\documentclass{amsart}

\usepackage{lscape,graphicx}
\usepackage{graphicx}

\usepackage{tikz}
\usetikzlibrary{arrows}
\tikzstyle{block}=[draw opacity=.7,line width=1.4cm]

\theoremstyle{definition}

\theoremstyle{remark}

\numberwithin{equation}{section}

\begin{document}

\title{Elegant expressions and generic formulas for the Riemann zeta function for integer arguments}

\author{{Michael A. Idowu}\\
\\
{SIMBIOS Centre, \\
School of Contemporary Sciences,\\
University of Abertay, \\
Dundee DD1 1HG, UK\\
m.idowu@abertay.ac.uk}}

\maketitle

\begin{abstract}
A new definition for the Riemann zeta function for positive integer number $s > 1$ is presented. We discover a most elegant expression and easy method for calculating the Riemann zeta function for small even integer values. Through this new reformulation we provide a one-line proof of $\zeta(2)={\pi^2\over{6}}$ and demonstrate that $\zeta(2s)$ may be calculated by hand using only the cotangent function when the magnitude of the integer s is small. 
\end{abstract}


\section{ Introduction}
More than a decade ago, K.S. K\"{o}lbig presented the polygamma function $\psi^{(k)}(z)$ for z=${1\over4}$ and z=${3\over4}$ \cite{Kol96} and remarked that he found it surprising that such values had not received much of the needed attention. Most intriguing is our latest discovery that beneath some simple operations involving the $\psi^{(k)}(z)$ function lies a fundamental architecture that firmly establishes a direct link between the Riemann zeta function and the well-known cotangent function. 

We present for the first time a new definition for $\zeta(s)$ for all positive integer number $s > 1$ and uncovers a most basic relationship between $\zeta(s)$ and the cotangent function $\cot(\pi z)$ at the point $z={1\over4}$ or $z={3\over4}$. We discover a most elegant expression and easy method for calculating the Riemann zeta function for small even integer values. Using this new method we show an original composition method of obtaining the actual values of $\zeta(2s)$ from $\cot(\pi z)$ and the constant $\pi$.  

\subsection{New definition and expression for $\zeta(s)$}

We manipulate the polygammma function $\psi^{(s)}(x)$ \cite{Nis11} and discover the relation 

\begin{equation}\label{main}
\zeta(s)= (-1)^s.{{( \psi^{(s-1)}({1\over4})+\psi^{(s-1)}({3\over4}))}\over{2^s.(2^s-1)}}.{1 \over{\Gamma(s)}}
\end{equation}
where s is a positive integer argument. This revelation shows that the polygamma function 

\begin{equation}
\psi^{(s-1)}(x) = {d^{s-1}\over{dx^{s-1}}}\psi(x) = {d^{s}\over{dx^{s}}}\ln \Gamma(x)
\end{equation}
may be used to calculate the Riemann zeta for any positive integer $s > 1$, e.g. 
\[
\zeta(2)= (-1)^2.{{( \psi^{'}({1\over4})+\psi^{'}({3\over4}))}\over{2^2.(2^2-1).\Gamma(2)}};
\zeta(3)= (-1)^3.{{( \psi^{''}({1\over4})+\psi^{''}({3\over4}))}\over{2^3.(2^3-1).\Gamma(3)}};
\]
\[
\zeta(4)= (-1)^4.{{( \psi^{'''}({1\over4})+\psi^{'''}({3\over4}))}\over{2^4.(2^4-1).\Gamma(4)}}; 
\zeta(5)= (-1)^5.{{( \psi^{''''}({1\over4})+\psi^{''''}({3\over4}))}\over{2^5.(2^5-1).\Gamma(5)}}; \dots 
\]

\section{Alternative generic formula for $\zeta(2s)$}
Given that the reflection relation of the polygamma function
\begin{equation}
(-1)^{s-1} \psi^{(s-1)}({1-z})-\psi^{(s-1)}({1-z}) = \pi{d^{(s-1)}\over{dz^{(s-1)}}}cot(\pi z) 
\end{equation}
holds, it is easy to infer from the relation  
\[
(-1)^{2s-1}\psi^{(2s-1)}({1-z})-\psi^{(2s-1)}({1-z}) = - (\psi^{(2s-1)}({1-z})+\psi^{(2s-1)}({1-z})) . 
\] that
\[
(-1)^{2s-1}\psi^{(2s-1)}({3\over4})-\psi^{(2s-1)}({1\over4}) = - (\psi^{(2s-1)}({3\over4})+\psi^{(2s-1)}({1\over4})), 
\]
 for all integer $s>0$.
Therefore
\begin{equation}
{2^{2s}(2^{2s}-1)\Gamma(2s)}\zeta(2s) = -\pi {d^{(2s-1)}\over{dz^{(2s-1)}}}cot(\pi z) \mid{_{z\rightarrow{1\over4}}}.
\end{equation}

\subsection{A one-line proof of $\zeta(2)={\pi^2\over6}$}
Calculating the values of zeta for any small even integer is easy and straightforward: for instance, we may derive the value of 
$\zeta(2)$ as follows:
\begin{equation}\label{temp} 
2^2(2^2-1)1!\zeta(2) = -\pi(-\pi(cot(\pi z)^2 + 1))\mid{_{z\rightarrow{1\over4}}} 
\end{equation}
\[ \Rightarrow \zeta(2) = {{(1+1)\pi^2}\over{3(4)}}= {{2\pi^2}\over{3(4)}} = {\pi^2\over6}.\]

\subsection{Derivation of the values of $\zeta(4)$, $\zeta(6)$, and $\zeta(8)$}
Other values of the zeta function for the subsequent even integers (e.g. $\zeta(4)$) may be calculated from taking twice derivatives of the last resultant cotangent function $-\pi(-\pi(cot(\pi z)^2 + 1))$ obtained in \ref{temp}, i.e. 
$-\pi{d\over{dz}}[{{d\over{dz}}{(-\pi(cot(\pi z)^2 + 1))}]}$:

\[ 2^4(2^4-1)3!\zeta(4) = -\pi ((-2)\pi^3(3cot(\pi z)^2 + 1)(cot(\pi z)^2 + 1)) \mid{_{z\rightarrow{1\over4}}} = (1+1)(3+1)2\pi^4; \]
\[\Rightarrow \zeta(4)={{2(2)4\pi^6}\over{16(15)6}} = {\pi^4\over{90}};\]

\[ 2^6(2^6-1)5!\zeta(6) = -\pi ((-8)\pi^5(cot(\pi z)^2 + 1)(15cot(\pi z)^4 + 15cot(\pi z)^2 + 2)) \mid{_{z\rightarrow{1\over4}}} \]
\[ \Rightarrow \zeta(6)={{(1+1)(15+15+2)8\pi^6}\over{64(63)120}} = {{8\pi^6}\over{(63)120}} = {\pi^6\over{945}};\]

\[ 2^8(2^8-1)5!\zeta(8) = -\pi ((-16)\pi^7(cot(\pi z)^2 + 1)(315 cot(\pi z)^6 + 525cot(\pi z)^4 + 231cot(\pi z)^2 + 17)) \mid{_{z\rightarrow{1\over4}}} \]
\[ \Rightarrow \zeta(8)={{(1+1)(315+525+231+17)16\pi^8}\over{256(255)5040}} = {{16\pi^8}\over{2(15)5040}} = {\pi^8\over{9450}};\]

\[\vdots\]
\begin{equation}\label{main2}
\zeta(2s) = -\pi {{d^{(2s-1)}\over{dz^{(2s-1)}}}cot(\pi z) \mid{_{z\rightarrow{1\over4}}}\over{{2^{2s}(2^{2s}-1)\Gamma(2s)}}}
\end{equation}

The formula presented in \ref{main2} may be used to derive the values of $\zeta(2s)$ by repeated differentiation, i.e. twice differentiating the resultant cotangent result from $\zeta(2s-2)$ to obtain the value of $\zeta(2s)$ as demonstrated. However, as s becomes larger this alternative method may become impracticable. In such cases, using the $\psi^{(s-1)}(z)$ method presented in \ref{main} might be more convenient.


\begin{thebibliography}{1}
\bibitem{Kol96}
K\"{o}lbig, K. S. ``The Polygamma Function $\psi^{(k)}(x)$ for $x={1\over{4}}$ and $x={3\over{4}}$." J. Comp. Appl. Math. 75, 43-46, 1996. 
\bibitem{Nis11}
Olver, Frank W. J., Lozier, Daniel W.,  Boisvert, Ronald F., Clark, Charles W. NIST Handbook of Mathematical Functions, Cambridge University Press (2010), ISBN 0521140633.
\end{thebibliography}
\end{document}